\documentclass[12pt,a4paper]{amsart}
\usepackage{amssymb,amsxtra}
\usepackage[english,french]{babel}

\theoremstyle{definition}

\theoremstyle{remark}

\theoremstyle{plain}

\theoremstyle{remark}

\newtheorem*{example}{Example}
\numberwithin{equation}{section}

\begin{document}

\title{On a conjecture regarding  enumeration of n-times persymmetric matrices over  $\mathbb{F}_{2} $ by rank }
\author{Jorgen~Cherly}
\address{D\'epartement de Math\'ematiques, Universit\'e de
    Brest, 29238 Brest cedex~3, France}
\email{Jorgen.Cherly@univ-brest.fr}
\email{andersen69@wanadoo.fr}

\maketitle

\begin{abstract}
Dans cet article nous annon\c{c}ons une conjecture concernant l'\' enum\' eration de n- fois matrices persym\' etriques sur  $\mathbb{F}_{2}$ par le rang.\\
Pour justifier notre assertion nous faisons remarquer que les formules obtenues sont valables pour n \' egal \`a un, deux et trois.
 \end{abstract}

\selectlanguage{english}

\begin{abstract}
In this paper we announce a conjecture concerning enumeration of  n-times  persymmetric  matrices over $\mathbb{F}_{2}$ by rank.
To justify our statement we remark that the formulas obtained are valid for n equal to one, two and three.
 \end{abstract}
\maketitle

\newpage
\tableofcontents
\newpage

\allowdisplaybreaks
\newpage

   \section{Some notations concerning the field of Laurent Series $ \mathbb{F}_{2}((T^{-1})) $ }
  \label{sec 1}  
  We denote by $ \mathbb{F}_{2}\big(\big({T^{-1}}\big) \big)
 = \mathbb{K} $ the completion
 of the field $\mathbb{F}_{2}(T), $  the field of  rational fonctions over the
 finite field\; $\mathbb{F}_{2}$,\; for the  infinity  valuation \;
 $ \mathfrak{v}=\mathfrak{v}_{\infty }$ \;defined by \;
 $ \mathfrak{v}\big(\frac{A}{B}\big) = degB -degA $ \;
 for each pair (A,B) of non-zero polynomials.
 Then every element non-zero t in
  $\mathbb{F}_{2}\big(\big({\frac{1}{T}}\big) \big) $
 can be expanded in a unique way in a convergent Laurent series
                              $  t = \sum_{j= -\infty }^{-\mathfrak{v}(t)}t_{j}T^j
                                 \; where\; t_{j}\in \mathbb{F}_{2}. $\\
  We associate to the infinity valuation\; $\mathfrak{v}= \mathfrak{v}_{\infty }$
   the absolute value \; $\vert \cdot \vert_{\infty} $\; defined by \;
  \begin{equation*}
  \vert t \vert_{\infty} =  \vert t \vert = 2^{-\mathfrak{v}(t)}. \\
\end{equation*}
    We denote  E the  Character of the additive locally compact group
$  \mathbb{F}_{2}\big(\big({\frac{1}{T}}\big) \big) $ defined by \\
\begin{equation*}
 E\big( \sum_{j= -\infty }^{-\mathfrak{v}(t)}t_{j}T^j\big)= \begin{cases}
 1 & \text{if      }   t_{-1}= 0, \\
  -1 & \text{if      }   t_{-1}= 1.
    \end{cases}
\end{equation*}
  We denote $\mathbb{P}$ the valuation ideal in $ \mathbb{K},$ also denoted the unit interval of  $\mathbb{K},$ i.e.
  the open ball of radius 1 about 0 or, alternatively, the set of all Laurent series 
   $$ \sum_{i\geq 1}\alpha _{i}T^{-i}\quad (\alpha _{i}\in  \mathbb{F}_{2} ) $$ and, for every rational
    integer j,  we denote by $\mathbb{P}_{j} $
     the  ideal $\left\{t \in \mathbb{K}|\; \mathfrak{v}(t) > j \right\}. $
     The sets\; $ \mathbb{P}_{j}$\; are compact subgroups  of the additive
     locally compact group \; $ \mathbb{K}. $\\
      All $ t \in \mathbb{F}_{2}\Big(\Big(\frac{1}{T}\Big)\Big) $ may be written in a unique way as
$ t = [t] + \left\{t\right\}, $ \;  $  [t] \in \mathbb{F}_{2}[T] ,
 \; \left\{t\right\}\in \mathbb{P}  ( =\mathbb{P}_{0}). $\\
 We denote by dt the Haar measure on \; $ \mathbb{K} $\; chosen so that \\
  $$ \int_{\mathbb{P}}dt = 1. $$\\
  
  $$ Let \quad
  (t_{1},t_{2},\ldots,t_{n} )
 =  \big( \sum_{j=-\infty}^{-\nu(t_{1})}\alpha _{j}^{(1)}T^{j},  \sum_{j=-\infty}^{-\nu(t_{2})}\alpha _{j}^{(2)}T^{j} ,\ldots, \sum_{j=-\infty}^{-\nu(t_{n})}\alpha _{j}^{(n)}T^{j}\big) \in  \mathbb{K}^{n}. $$ 
 We denote $\psi  $  the  Character on  $(\mathbb{K}^n, +) $ defined by \\
 \begin{align*}
  \psi \big( \sum_{j=-\infty}^{-\nu(t_{1})}\alpha _{j}^{(1)}T^{j},  \sum_{j=-\infty}^{-\nu(t_{2})}\alpha _{j}^{(2)}T^{j} ,\ldots, \sum_{j=-\infty}^{-\nu(t_{n})}\alpha _{j}^{(n)}T^{j}\big) & = E \big( \sum_{j=-\infty}^{-\nu(t_{1})}\alpha _{j}^{(1)}T^{j}\big) \cdot E\big( \sum_{j=-\infty}^{-\nu(t_{2})}\alpha _{j}^{(2)}T^{j}\big)\cdots E\big(  \sum_{j=-\infty}^{-\nu(t_{n})}\alpha _{j}^{(n)}T^{j}\big) \\
  & = 
    \begin{cases}
 1 & \text{if      }     \alpha _{-1}^{(1)} +    \alpha _{-1}^{(2)}  + \ldots +   \alpha _{-1}^{(n)}   = 0 \\
  -1 & \text{if      }    \alpha _{-1}^{(1)} +    \alpha _{-1}^{(2)}  + \ldots +   \alpha _{-1}^{(n)}   =1                                                                                                                          
    \end{cases}
  \end{align*}
        \subsection{Computation of the number  $ \Gamma_{i}^{\left[1\atop{\vdots \atop 1}\right]\times 2} $   of   n-times persymmetric 
   $n \times 2$  rank i matrices }
  \label{subsec 1} 
    $$ Set\quad
  (t_{1},t_{2},\ldots,t_{n} )
 =  \big( \sum_{i\geq 1}\alpha _{i}^{(1)}T^{-i}, \sum_{i \geq 1}\alpha  _{i}^{(2)}T^{-i},\sum_{i \geq 1}\alpha _{i}^{(3)}T^{-i},\ldots,\sum_{i \geq 1}\alpha _{i}^{(n)}T^{-i}   \big) \in  \mathbb{P}^{n}. $$ 
 We denote by $D^{\left[1\atop{\vdots \atop 1}\right]\times 2}(t_{1},t_{2},\ldots,t_{n} ) $
    the following $n \times 2$ \;  n-times  persymmetric  matrix  over the finite field  $\mathbb{F}_{2} $ 
       \begin{displaymath}
   \left (  \begin{array} {cc}
\alpha  _{1}^{(1)} & \alpha  _{2}^{(1)}   \\
\hline \\
\alpha  _{1}^{(2)} & \alpha  _{2}^{(2)}   \\
\hline\\
\alpha  _{1}^{(3)} & \alpha  _{2}^{(3)}  \\
\hline \\
\vdots & \vdots \\
\hline \\
\alpha  _{1}^{(n)} & \alpha  _{2}^{(n)}  \\
\end{array} \right )   \; \overset{\text{rank}}{\sim}  \;
 \left (  \begin{array} {cccc}
     \alpha  _{1} & \alpha  _{2}  & \ldots  &  \alpha _{n} \\
          \beta _{1} &  \beta_{2}  & \ldots  &   \beta_{n} 
  \end{array} \right )
\end{displaymath} 
  Let $ \displaystyle  f (t_{1},t_{2},\ldots,t_{n} ) $  be the exponential sum  in $ \mathbb{P}^{n} $ defined by\\
    $(t_{1},t_{2},\ldots,t_{n} ) \displaystyle\in \mathbb{P}^{n}\longrightarrow \\
    \sum_{deg Y\leq 1}\sum_{deg U_{1}\leq  0}E(t_{1} YU_{1})
  \sum_{deg U_{2} \leq 0}E(t _{2} YU_{2}) \ldots \sum_{deg U_{n} \leq 0} E(t _{n} YU_{n}). $\vspace{0.5 cm}\\
  Then
  $$     f (t_{1},t_{2},\ldots,t_{n} ) =
  2^{n+2- rank\big[ D^{\left[1\atop{\vdots \atop 1}\right]\times 2}(t_{1},t_{2},\ldots,t_{n} )  \big ] } $$

  Hence  the number denoted by $ R_{q} $ of solutions \\
  
 $(Y_1,U_{1}^{(1)},U_{2}^{(1)}, \ldots,U_{n}^{(1)}, Y_2,U_{1}^{(2)},U_{2}^{(2)}, 
\ldots,U_{n}^{(2)},\ldots  Y_q,U_{1}^{(q)},U_{2}^{(q)}, \ldots,U_{n}^{(q)}   ) \in ( \mathbb{F}_{2}[T])^{(n+1)q} $ \vspace{0.5 cm}\\
 of the polynomial equations  \vspace{0.5 cm}
  \[\left\{\begin{array}{c}
 Y_{1}U_{1}^{(1)} + Y_{2}U_{1}^{(2)} + \ldots  + Y_{q}U_{1}^{(q)} = 0  \\
    Y_{1}U_{2}^{(1)} + Y_{2}U_{2}^{(2)} + \ldots  + Y_{q}U_{2}^{(q)} = 0\\
    \vdots \\
   Y_{1}U_{n}^{(1)} + Y_{2}U_{n}^{(2)} + \ldots  + Y_{q}U_{n}^{(q)} = 0 
 \end{array}\right.\]
    satisfying the degree conditions \\
                   $$  degY_i \leq 1 ,
                   \quad degU_{j}^{(i)} \leq 0, \quad  for \quad 1\leq j\leq n  \quad 1\leq i \leq q $$ \\
  is equal to the following integral over the unit interval in $ \mathbb{K}^{n} $
    $$ \int_{\mathbb{P}^{n}} f^{q}(t_{1},t_{2},\ldots,t_{n}) dt_{1}dt _{2}\ldots dt _{n}. $$
  Observing that $ f (t_{1},t_{2},\ldots,t_{n} ) $ is constant on cosets of $ \mathbb{P}_{2}^{n} $ in $ \mathbb{P}^{n} $\;
  the above integral is equal to 
  
  \begin{equation}
  \label{eq 1.1}
 2^{q(n+2) - 2n}\sum_{i = 0}^{2}
 \Gamma_{i}^{\left[1\atop{\vdots \atop 1}\right]\times 2} 2^{-iq} =  R_{q} 
 \end{equation}
 
 From \eqref{eq 1.1} we obtain for q = 1\\
 
   \begin{equation}
  \label{eq 1.2}
 2^{2-n}\sum_{i = 0}^{2}
 \Gamma_{i}^{\left[1\atop{\vdots \atop 1}\right]\times 2} 2^{-i} = 2^{n}+2^2-1
 \end{equation}

We have obviously \\

   \begin{equation}
  \label{eq 1.3}
 \sum_{i = 0}^{2}
 \Gamma_{i}^{\left[1\atop{\vdots \atop 1}\right]\times 2}  = 2^{2n}
 \end{equation}

  Combining \eqref{eq 1.2}, \eqref{eq 1.3} we get \\
  
         \begin{equation}
         \label{eq 1.4}
 \Gamma_{i}^{\left[1\atop{\vdots \atop 1}\right]\times 2} =  \begin{cases}
1 & \text{if  } i = 0,        \\
 (2^{n}-1)\cdot 3 & \text{if   } i=1,\\
 2^{2n} -3\cdot2^{n}+2 & \text{if   }  i = 2
    \end{cases}    
   \end{equation}
   
From     \eqref{eq 1.1}, \eqref{eq 1.4} we obtain :\\

\begin{equation*}
R_{q} = 2^{(q-2)n}\cdot [2^{2q} +2^{2n} +3\cdot(2^{n+q}-2^{q}-2^{n}) +2]
\end{equation*}

  \subsection{Computation of the number  $ \Gamma_{i}^{\left[2\atop{\vdots \atop 2}\right]\times 3} $   of   n-times persymmetric 
   $2n \times 3$  rank i matrices }
  \label{subsec 2}  
 $$ Set\quad
  (t_{1},t_{2},\ldots,t_{n} )
 =  \big( \sum_{i\geq 1}\alpha _{i}^{(1)}T^{-i}, \sum_{i \geq 1}\alpha  _{i}^{(2)}T^{-i},\sum_{i \geq 1}\alpha _{i}^{(3)}T^{-i},\ldots,\sum_{i \geq 1}\alpha _{i}^{(n)}T^{-i}   \big) \in  \mathbb{P}^{n}. $$ 
  We denote by $D^{\left[2\atop{\vdots \atop 2}\right]\times 3}(t_{1},t_{2},\ldots,t_{n} ) $
    the following $2n \times 3 $ \;  n-times  persymmetric  matrix  over the finite field  $\mathbb{F}_{2} $ 
   \begin{displaymath}
   \left (  \begin{array} {ccc}
\alpha  _{1}^{(1)} & \alpha  _{2}^{(1)}  &   \alpha_{3}^{(1)}   \\
\alpha  _{2}^{(1)} & \alpha  _{3}^{(1)}  &   \alpha_{4}^{(1)}  \\ 
\hline \\
\alpha  _{1}^{(2)} & \alpha  _{2}^{(2)}  &   \alpha_{3}^{(2)}   \\
\alpha  _{2}^{(2)} & \alpha  _{3}^{(2)}  &   \alpha_{4}^{(2)}  \\ 
\hline\\
\alpha  _{1}^{(3)} & \alpha  _{2}^{(3)}  &   \alpha_{3}^{(3)}   \\
\alpha  _{2}^{(3)} & \alpha  _{3}^{(3)}  &   \alpha_{4}^{(3)}  \\ 
\hline \\
\vdots & \vdots & \vdots \\
\hline \\
\alpha  _{1}^{(n)} & \alpha  _{2}^{(n)}  &   \alpha_{3}^{(n)}   \\
\alpha  _{2}^{(n)} & \alpha  _{3}^{(n)}  &   \alpha_{4}^{(n)}  \\ 
\end{array} \right )  
\end{displaymath} 

  Let $ \displaystyle  f (t_{1},t_{2},\ldots,t_{n} ) $  be the exponential sum  in $ \mathbb{P}^{n} $ defined by\\
    $(t_{1},t_{2},\ldots,t_{n} ) \displaystyle\in \mathbb{P}^{n}\longrightarrow \\
    \sum_{deg Y\leq 2}\sum_{deg U_{1}\leq  1}E(t_{1} YU_{1})
  \sum_{deg U_{2} \leq 1}E(t _{2} YU_{2}) \ldots \sum_{deg U_{n} \leq 1} E(t _{n} YU_{n}). $\vspace{0.5 cm}\\
  Then
  $$     f (t_{1},t_{2},\ldots,t_{n} ) =
  2^{2n+3- rank\big[ D^{\left[2\atop{\vdots \atop 2}\right]\times 3}(t_{1},t_{2},\ldots,t_{n} )  \big ] } $$

  Hence  the number denoted by $ R_{q} $ of solutions \\
  
 $(Y_1,U_{1}^{(1)},U_{2}^{(1)}, \ldots,U_{n}^{(1)}, Y_2,U_{1}^{(2)},U_{2}^{(2)}, 
\ldots,U_{n}^{(2)},\ldots  Y_q,U_{1}^{(q)},U_{2}^{(q)}, \ldots,U_{n}^{(q)}   ) \in ( \mathbb{F}_{2}[T])^{(n+1)q} $ \vspace{0.5 cm}\\
 of the polynomial equations  \vspace{0.5 cm}
  \[\left\{\begin{array}{c}
 Y_{1}U_{1}^{(1)} + Y_{2}U_{1}^{(2)} + \ldots  + Y_{q}U_{1}^{(q)} = 0  \\
    Y_{1}U_{2}^{(1)} + Y_{2}U_{2}^{(2)} + \ldots  + Y_{q}U_{2}^{(q)} = 0\\
    \vdots \\
   Y_{1}U_{n}^{(1)} + Y_{2}U_{n}^{(2)} + \ldots  + Y_{q}U_{n}^{(q)} = 0 
 \end{array}\right.\]
    satisfying the degree conditions \\
                   $$  degY_i \leq 2 ,
                   \quad degU_{j}^{(i)} \leq 1, \quad  for \quad 1\leq j\leq n  \quad 1\leq i \leq q $$ \\
  is equal to the following integral over the unit interval in $ \mathbb{K}^{n} $
    $$ \int_{\mathbb{P}^{n}} f^{q}(t_{1},t_{2},\ldots,t_{n}) dt_{1}dt _{2}\ldots dt _{n}. $$
  Observing that $ f (t_{1},t_{2},\ldots,t_{n} ) $ is constant on cosets of $ \mathbb{P}_{4}^{n} $ in $ \mathbb{P}^{n} $\;
  the above integral is equal to 
  
  \begin{equation}
  \label{eq 1.5}
 2^{q(2n+3) - 4n}\sum_{i = 0}^{3}
 \Gamma_{i}^{\left[2\atop{\vdots \atop 2}\right]\times 3} 2^{-iq} =  R_{q} 
 \end{equation}
 
 From \eqref{eq 1.5} we obtain for q = 1\\
 
   \begin{equation}
  \label{eq 1.6}
 2^{3-2n}\sum_{i = 0}^{3}
 \Gamma_{i}^{\left[2\atop{\vdots \atop 2}\right]\times 3} 2^{-i} = 2^{2n}+2^3-1
 \end{equation}

We have obviously \\

   \begin{equation}
  \label{eq 1.7}
 \sum_{i = 0}^{3}
 \Gamma_{i}^{\left[2\atop{\vdots \atop 2}\right]\times 3}  = 2^{4n}
 \end{equation}

From  the fact that the number of rank one persymmetric  matrices over $\mathbb{F}_{2}$ is equal to three  we obtain using
 combinatorial methods  that: \\
 
    \begin{equation}
  \label{eq 1.8}
 \Gamma_{1}^{\left[2\atop{\vdots \atop 2}\right]\times 3}  = (2^{n}-1)\cdot 3
 \end{equation}

  Combining \eqref{eq 1.6}, \eqref{eq 1.7} and \eqref{eq 1.8} we get \\
  
         \begin{equation}
         \label{eq 1.9}
 \Gamma_{i}^{\left[2\atop{\vdots \atop 2}\right]\times 3} =  \begin{cases}
1 & \text{if  } i = 0,        \\
 (2^{n}-1)\cdot 3 & \text{if   } i=1,\\
 7\cdot 2^{2n} -9\cdot2^{n}+2 & \text{if   }  i = 2,  \\
2^{4n} - 7\cdot 2^{2n} +6\cdot2^{n}  & \text{if   }  i=3 
    \end{cases}    
   \end{equation}

  \textbf{Generalization:}
  
    Let $s_{j} \geqslant 2 $ for $1\leqslant j \leqslant n ,$ denote by $D^{\left[s_{1}\atop{\vdots \atop s_{n}}\right]\times 3}(t_{1},t_{2},\ldots,t_{n} ) $
    
    the following $(\sum_{j=1}^{n}s_{j}) \times 3 $ \;  n-times  persymmetric  matrix  over the finite field  $\mathbb{F}_{2} $ 
    
  \begin{displaymath}
   \left (  \begin{array} {ccc}
\alpha  _{1}^{(1)} & \alpha  _{2}^{(1)}  &   \alpha_{3}^{(1)}   \\
\alpha  _{2}^{(1)} & \alpha  _{3}^{(1)}  &   \alpha_{4}^{(1)}  \\ 
\vdots & \vdots & \vdots \\
\alpha  _{s_{1}}^{(1)} & \alpha  _{s_{1}+1}^{(1)}  &   \alpha_{s_{1}+2}^{(1)}  \\ 
\hline \\
\alpha  _{1}^{(2)} & \alpha  _{2}^{(2)}  &   \alpha_{3}^{(2)}   \\
\alpha  _{2}^{(2)} & \alpha  _{3}^{(2)}  &   \alpha_{4}^{(2)}  \\ 
\vdots & \vdots & \vdots \\
\alpha  _{s_{2}}^{(2)} & \alpha  _{s_{2}+1}^{(2)}  &   \alpha_{s_{2}+2}^{(2)}  \\ 
\hline\\
\alpha  _{1}^{(3)} & \alpha  _{2}^{(3)}  &   \alpha_{3}^{(3)}   \\
\alpha  _{2}^{(3)} & \alpha  _{3}^{(3)}  &   \alpha_{4}^{(3)}  \\ 
\vdots & \vdots & \vdots \\
\alpha  _{s_{3}}^{(3)} & \alpha  _{s_{3}+1}^{(3)}  &   \alpha_{s_{3}+2}^{(3)}  \\ 
\hline \\
\vdots & \vdots & \vdots \\
\hline \\
\alpha  _{1}^{(n)} & \alpha  _{2}^{(n)}  &   \alpha_{3}^{(n)}   \\
\alpha  _{2}^{(n)} & \alpha  _{3}^{(n)}  &   \alpha_{4}^{(n)}  \\ 
\vdots & \vdots & \vdots \\
\alpha  _{s_{n}}^{(n)} & \alpha  _{s_{n}+1}^{(n)}  &   \alpha_{s_{n}+2}^{(n)}  \\ 
\end{array} \right )  
\end{displaymath}

  Let $ \displaystyle  f (t_{1},t_{2},\ldots,t_{n} ) $  be the exponential sum  in $ \mathbb{P}^{n} $ defined by\\
    $(t_{1},t_{2},\ldots,t_{n} ) \displaystyle\in \mathbb{P}^{n}\longrightarrow \\
    \sum_{deg Y\leq 2}\sum_{deg U_{1}\leq  s_{1}-1}E(t_{1} YU_{1})
  \sum_{deg U_{2} \leq s_{2}-1}E(t _{2} YU_{2}) \ldots \sum_{deg U_{n} \leq s_{n}-1} E(t _{n} YU_{n}). $\vspace{0.5 cm}\\
    Then
  $$     f (t_{1},t_{2},\ldots,t_{n} ) =
  2^{\sum_{i=1}^{n}s_{i}+3- rank\big[ D^{\left[s_{1}\atop{\vdots \atop s_{n}}\right]\times 3}(t_{1},t_{2},\ldots,t_{n} )\big] } $$
    Hence  the number denoted by $ R_{q} $ of solutions \\
  
 $(Y_1,U_{1}^{(1)},U_{2}^{(1)}, \ldots,U_{n}^{(1)}, Y_2,U_{1}^{(2)},U_{2}^{(2)}, 
\ldots,U_{n}^{(2)},\ldots  Y_q,U_{1}^{(q)},U_{2}^{(q)}, \ldots,U_{n}^{(q)}   ) $ \vspace{0.5 cm}\\
 of the polynomial equations  \vspace{0.5 cm}
  \[\left\{\begin{array}{c}
 Y_{1}U_{1}^{(1)} + Y_{2}U_{1}^{(2)} + \ldots  + Y_{q}U_{1}^{(q)} = 0  \\
    Y_{1}U_{2}^{(1)} + Y_{2}U_{2}^{(2)} + \ldots  + Y_{q}U_{2}^{(q)} = 0\\
    \vdots \\
   Y_{1}U_{n}^{(1)} + Y_{2}U_{n}^{(2)} + \ldots  + Y_{q}U_{n}^{(q)} = 0 
 \end{array}\right.\]
    satisfying the degree conditions \\
                   $$  degY_i \leq 2 ,
                   \quad degU_{j}^{(i)} \leq s_{j}-1, \quad  for \quad 1\leq j\leq n  \quad 1\leq i \leq q $$ \\
  is equal to the following integral over the unit interval in $ \mathbb{K}^{n} $
    $$ \int_{\mathbb{P}^{n}} f^{q}(t_{1},t_{2},\ldots,t_{n}) dt_{1}dt _{2}\ldots dt _{n}. $$
  Observing that $ f (t_{1},t_{2},\ldots,t_{n} ) $ is constant on cosets of $ \prod_{j=1}^{n}\mathbb{P}_{s_{j}+2} $ in $ \mathbb{P}^{n} $\;
  the above integral is equal to 
  
  \begin{equation}
  \label{eq 1.10}
 2^{q(\sum_{j=1}^{n}s_{j}+3) -\sum_{j=1}^{n}s_{j}-2\cdot n}\sum_{i = 0}^{3}
  \Gamma_{i}^{\left[s_{1}\atop{\vdots \atop s_{n}}\right]\times 3} 2^{-iq} =  R_{q} 
 \end{equation}

 From \eqref{eq 1.10} we obtain for q = 1\\
 
   \begin{equation}
  \label{eq 1.11}
 2^{3-2n}\sum_{i = 0}^{3}
 \Gamma_{i}^{\left[s_{1}\atop{\vdots \atop s_{n}}\right]\times 3} 2^{-i} = 2^{\sum_{j=1}^{n}s_{j}}+2^3-1
 \end{equation}

We have obviously \\

   \begin{equation}
  \label{eq 1.12}
 \sum_{i = 0}^{3}
 \Gamma_{i}^{\left[s_{1}\atop{\vdots \atop s_{1}}\right]\times 3}  = 2^{\sum_{j=1}^{n}s_{j}+2n}
 \end{equation}

From  the fact that the number of rank one persymmetric  matrices over $\mathbb{F}_{2}$ is equal to three  we obtain as above using
 combinatorial methods  : \\
 
    \begin{equation}
  \label{eq 1.13}
 \Gamma_{1}^{\left[s_{1}\atop{\vdots \atop s_{n}}\right]\times 3}  = (2^{n}-1)\cdot 3
 \end{equation}

  Combining \eqref{eq 1.11}, \eqref{eq 1.12} and \eqref{eq 1.13} we get \\
  
         \begin{equation}
         \label{eq 1.14}
 \Gamma_{i}^{\left[s_{1}\atop{\vdots \atop s_{n}}\right]\times 3} =  \begin{cases}
1 & \text{if  } i = 0,        \\
 (2^{n}-1)\cdot 3 & \text{if   } i=1,\\
 7\cdot 2^{2n} -9\cdot2^{n}+2 & \text{if   }  i = 2,  \\
2^{  \sum_{j=1}^{n}s_{j} + 2n} - 7\cdot 2^{2n} +6\cdot2^{n}  & \text{if   }  i=3 
    \end{cases}    
   \end{equation}

We get  from  \eqref{eq 1.14}, \eqref{eq 1.9} and \eqref{eq 1.4}  whenever  $s_{j} \geqslant 2 $ for  $1\leqslant j \leqslant n $ \\

     \begin{equation}
         \label{eq 1.15}
 \Gamma_{i}^{\left[s_{1}\atop{\vdots \atop s_{n}}\right]\times 3} =  \begin{cases}
1 & \text{if  } i = 0,        \\
 \Gamma_{1}^{\left[1\atop{\vdots \atop 1}\right]\times 2}  =   (2^{n}-1)\cdot 3 & \text{if   } i=1,\\
  \Gamma_{2}^{\left[2\atop{\vdots \atop 2}\right]\times 3} =    7\cdot 2^{2n} -9\cdot2^{n}+2 & \text{if   }  i = 2,  \\
2^{  \sum_{j=1}^{n}s_{j} + 2n} - 7\cdot 2^{2n} +6\cdot2^{n}  & \text{if   }  i=3 
    \end{cases}    
   \end{equation}
  \begin{example}
    
   We obtain from \eqref{eq 1.9} with n=4:\\
      \begin{equation*}
  \Gamma_{i}^{\left[2\atop {2\atop {2\atop2 }}\right]\times 3} =  \begin{cases}
1 & \text{if  } i = 0,        \\
 45 & \text{if   } i=1,\\
 1650 & \text{if   }  i = 2,  \\
  63840 & \text{if   }  i=3 
    \end{cases}    
   \end{equation*}
 The number $ \Gamma_{i}^{\left[2\atop {2\atop {2\atop{2 \atop (1)}}}\right]\times 3} $ of rank i matrices of the form \\
   \begin{displaymath}
 \left (  \begin{array} {ccccccccc}
     \alpha  _{1} & \alpha  _{2}  &   \alpha_{3} &  \alpha _{4}  &  \alpha  _{5} & \alpha  _{6}  &   \alpha_{7} &  \alpha _{8}  &  \alpha _{9} \\
        \alpha  _{6} & \alpha  _{7}  &   \alpha_{8} &  \alpha _{19}  &  \alpha  _{10} & \alpha  _{11}  &   \alpha_{12} &  \alpha _{13} &  \alpha _{14}   \\
           \alpha  _{11} & \alpha  _{12}  &   \alpha_{13} &  \alpha _{14}  &  \alpha  _{15} & \alpha  _{16}  &   \alpha_{17} &  \alpha _{18} & \alpha_{19}    
                  \end{array} \right ) \; \overset{\text{rank}}{\sim}  \;
 \left (  \begin{array} {ccc}
\alpha  _{1} & \alpha  _{2}  &   \alpha_{3}   \\
\alpha  _{2} & \alpha  _{3}  &   \alpha_{4}  \\ 
\hline \\
\beta  _{1} & \beta _{2}  &   \beta_{3}   \\
\beta  _{2} & \beta  _{3}  &   \beta_{4}  \\ 
\hline \\
\gamma _{1} & \gamma _{2} & \gamma  _{3}\\
 \gamma _{2} & \gamma  _{3} & \gamma _{4} \\
\hline \\
\mu _{1} & \mu _{2} & \mu  _{3}\\
 \mu_{2} & \mu _{3} & \mu _{4} \\
\hline \\
\delta_{11} & \delta_{12} & \delta  _{13}
\end{array} \right )    
\end{displaymath} \\
is equal to \vspace{0.1 cm} \\

  $ 2^{i}\cdot  \Gamma_{i}^{\left[2\atop {2\atop {2\atop2 }}\right]\times 3}
  +  (2^{3}-2^{i-1})\cdot \Gamma_{i}^{\left[2\atop {2\atop {2\atop2 }}\right]\times 3}  for\quad 0\leq i\leq inf(3,9)  \quad \text{see} \; [2]$  
  \vspace{0.1 cm} \\
  That is \\

      \begin{equation*}
\Gamma_{i}^{\left[2\atop {2\atop {2\atop{2 \atop (1)}}}\right]\times 3}  =  \begin{cases}
1 & \text{if  } i = 0,        \\
 97 & \text{if   } i=1,\\
 6870 & \text{if   }  i = 2,  \\
  5177320 & \text{if   }  i=3 
    \end{cases}    
   \end{equation*}
 \end{example}

   \subsection{Computation of the number  $ \Gamma_{i}^{\left[s_{1}\atop{\vdots \atop s_{n}}\right]\times 4} $   of   n-times persymmetric 
   $(\sum_{j=1}^{n}s_{j}) \times 4$  rank i matrices  whenever  $s_{j} \geqslant 3  $ for  $1\leqslant j \leqslant n $ }
  \label{subsec 3}

  Denote by $D^{\left[s_{1}\atop{\vdots \atop s_{n}}\right]\times 4}(t_{1},t_{2},\ldots,t_{n} ) $
    
    the following $(\sum_{j=1}^{n}s_{j}) \times 4 $ \;  n-times  persymmetric  matrix  over the finite field  $\mathbb{F}_{2} $ 
    
  \begin{displaymath}
   \left (  \begin{array} {cccc}
\alpha  _{1}^{(1)} & \alpha  _{2}^{(1)}  &   \alpha_{3}^{(1)} &   \alpha_{4}^{(1)}  \\
\alpha  _{2}^{(1)} & \alpha  _{3}^{(1)}  &   \alpha_{4}^{(1)} &   \alpha_{5}^{(1)} \\ 
\vdots & \vdots & \vdots \\
\alpha  _{s_{1}}^{(1)} & \alpha  _{s_{1}+1}^{(1)}  &   \alpha_{s_{1}+2}^{(1)} &   \alpha_{s_{1}+3}^{(1)}  \\ 
\hline \\
\alpha  _{1}^{(2)} & \alpha  _{2}^{(2)}  &   \alpha_{3}^{(2)} &   \alpha_{4}^{(2)}  \\
\alpha  _{2}^{(2)} & \alpha  _{3}^{(2)}  &   \alpha_{4}^{(2)} &   \alpha_{5}^{(2)} \\ 
\vdots & \vdots & \vdots \\
\alpha  _{s_{2}}^{(2)} & \alpha  _{s_{2}+1}^{(2)}  &   \alpha_{s_{2}+2}^{(2)} &   \alpha_{s_{2}+3}^{(2)}  \\ 
\hline\\
\alpha  _{1}^{(3)} & \alpha  _{2}^{(3)}  &   \alpha_{3}^{(3)}  &   \alpha_{4}^{(3)} \\
\alpha  _{2}^{(3)} & \alpha  _{3}^{(3)}  &   \alpha_{4}^{(3)}&   \alpha_{5}^{(3)}  \\ 
\vdots & \vdots & \vdots \\
\alpha  _{s_{3}}^{(3)} & \alpha  _{s_{3}+1}^{(3)}  &   \alpha_{s_{3}+2}^{(3)}&   \alpha_{s_{3}+3}^{(3)}   \\ 
\hline \\
\vdots & \vdots & \vdots \\
\hline \\
\alpha  _{1}^{(n)} & \alpha  _{2}^{(n)}  &   \alpha_{3}^{(n)} &   \alpha_{4}^{(n)}  \\
\alpha  _{2}^{(n)} & \alpha  _{3}^{(n)}  &   \alpha_{4}^{(n)}&   \alpha_{5}^{(n)}  \\ 
\vdots & \vdots & \vdots \\
\alpha  _{s_{n}}^{(n)} & \alpha  _{s_{n}+1}^{(n)}  &   \alpha_{s_{n}+2}^{(n)} &   \alpha_{s_{n}+3}^{(n)}  \\ 
\end{array} \right )  
\end{displaymath}

  Let $ \displaystyle  f (t_{1},t_{2},\ldots,t_{n} ) $  be the exponential sum  in $ \mathbb{P}^{n} $ defined by\\
    $(t_{1},t_{2},\ldots,t_{n} ) \displaystyle\in \mathbb{P}^{n}\longrightarrow \\
    \sum_{deg Y\leq 3}\sum_{deg U_{1}\leq  s_{1}-1}E(t_{1} YU_{1})
  \sum_{deg U_{2} \leq s_{2}-1}E(t _{2} YU_{2}) \ldots \sum_{deg U_{n} \leq s_{n}-1} E(t _{n} YU_{n}). $\vspace{0.5 cm}\\
    Then
  $$     f (t_{1},t_{2},\ldots,t_{n} ) =
  2^{\sum_{j=1}^{n}s_{j}+4- rank\big[ D^{\left[s_{1}\atop{\vdots \atop s_{n}}\right]\times 4}(t_{1},t_{2},\ldots,t_{n} )\big] } $$
    Hence  the number denoted by $ R_{q} $ of solutions \\
  
 $(Y_1,U_{1}^{(1)},U_{2}^{(1)}, \ldots,U_{n}^{(1)}, Y_2,U_{1}^{(2)},U_{2}^{(2)}, 
\ldots,U_{n}^{(2)},\ldots  Y_q,U_{1}^{(q)},U_{2}^{(q)}, \ldots,U_{n}^{(q)}   ) $ \vspace{0.5 cm}\\
 of the polynomial equations  \vspace{0.5 cm}
  \[\left\{\begin{array}{c}
 Y_{1}U_{1}^{(1)} + Y_{2}U_{1}^{(2)} + \ldots  + Y_{q}U_{1}^{(q)} = 0  \\
    Y_{1}U_{2}^{(1)} + Y_{2}U_{2}^{(2)} + \ldots  + Y_{q}U_{2}^{(q)} = 0\\
    \vdots \\
   Y_{1}U_{n}^{(1)} + Y_{2}U_{n}^{(2)} + \ldots  + Y_{q}U_{n}^{(q)} = 0 
 \end{array}\right.\]
    satisfying the degree conditions \\
                   $$  degY_i \leq 3 ,
                   \quad degU_{j}^{(i)} \leq s_{j}-1, \quad  for \quad 1\leq j\leq n  \quad 1\leq i \leq q $$ \\
  is equal to the following integral over the unit interval in $ \mathbb{K}^{n} $
    $$ \int_{\mathbb{P}^{n}} f^{q}(t_{1},t_{2},\ldots,t_{n}) dt_{1}dt _{2}\ldots dt _{n}. $$
  Observing that $ f (t_{1},t_{2},\ldots,t_{n} ) $ is constant on cosets of $ \prod_{j=1}^{n}\mathbb{P}_{s_{j}+3} $ in $ \mathbb{P}^{n} $\;
  the above integral is equal to 
  
  \begin{equation}
  \label{eq 1.16}
 2^{q(\sum_{j=1}^{n}s_{j}+4) -\sum_{j=1}^{n}s_{j}-3\cdot n}\sum_{i = 0}^{4}
  \Gamma_{i}^{\left[s_{1}\atop{\vdots \atop s_{n}}\right]\times 4} 2^{-iq} =  R_{q} 
 \end{equation}

 From \eqref{eq 1.16} we obtain for q = 1\\
 
   \begin{equation}
  \label{eq 1.17}
 2^{4-3n}\sum_{i = 0}^{4}
 \Gamma_{i}^{\left[s_{1}\atop{\vdots \atop s_{n}}\right]\times 4} 2^{-i} = 2^{\sum_{j=1}^{n}s_{j}}+2^4-1
 \end{equation}

We have obviously \\

   \begin{equation}
  \label{eq 1.18}
 \sum_{i = 0}^{4}
 \Gamma_{i}^{\left[s_{1}\atop{\vdots \atop s_{1}}\right]\times 4}  = 2^{\sum_{j=1}^{n}s_{j}+3n}
 \end{equation}

From  the fact that the number of rank one persymmetric  matrices over $\mathbb{F}_{2}$ is equal to three  we obtain as above using
 combinatorial methods  : \\
 
    \begin{equation}
  \label{eq 1.19}
 \Gamma_{1}^{\left[s_{1}\atop{\vdots \atop s_{n}}\right]\times 4}  = (2^{n}-1)\cdot 3
 \end{equation}

 We assume now: \\
 
    \begin{equation}
         \label{eq 1.20}
 \Gamma_{i}^{\left[s_{1}\atop{\vdots \atop s_{n}}\right]\times 4} =  \begin{cases}
1 & \text{if  } i = 0,        \\
 \Gamma_{1}^{\left[1\atop{\vdots \atop 1}\right]\times 2}  =   (2^{n}-1)\cdot 3 & \text{if   } i=1,\\
  \Gamma_{2}^{\left[2\atop{\vdots \atop 2}\right]\times 3} =    7\cdot 2^{2n} -9\cdot2^{n}+2 & \text{if   }  i = 2,  \\
 \end{cases}    
   \end{equation}
To justify our assumption \eqref{eq 1.20} for i=2 we remark that:\\
\begin{itemize}
\item  The number  of rank two persymmetric  matrices over $\mathbb{F}_{2}$ is equal to $ 7\cdot 2^{2n} -9\cdot2^{n}+2 = 7\cdot 2^{2} -9\cdot2^{1}+2 = 12  $   for n=1 [see (1), (2)]  \\
\item The number  of rank two double  persymmetric  matrices over $\mathbb{F}_{2}$ is equal to $ 7\cdot 2^{2n} -9\cdot2^{n}+2 = 7\cdot 2^{4} -9\cdot2^{2}+2 = 78  $   for n=2 [see (3)]  \\
\item The number  of rank two  triple   persymmetric  matrices over $\mathbb{F}_{2}$ is equal to $ 7\cdot 2^{2n} -9\cdot2^{n}+2 = 7\cdot 2^{6} -9\cdot2^{3}+2 = 378 $   for n=3 [see (4)]  \\
\end{itemize}

  Combining \eqref{eq 1.18}, \eqref{eq 1.19} and \eqref{eq 1.20} we state that the number $ \Gamma_{i}^{\left[s_{1}\atop{\vdots \atop s_{n}}\right]\times 4} $   of   n-times persymmetric 
   $(\sum_{j=1}^{n}s_{j}) \times 4$  rank i matrices  is equal to : \\ 
  
         \begin{equation}
         \label{eq 1.21}
  \begin{cases}
1 & \text{if  } i = 0,        \\
 (2^{n}-1)\cdot 3 & \text{if   } i=1,\\
 7\cdot 2^{2n} -9\cdot2^{n}+2 & \text{if   }  i = 2,  \\
  15\cdot 2^{3n} -21\cdot2^{2n}+3\cdot2^{n+1} & \text{if   }  i = 3,  \\
2^{  \sum_{j=1}^{n}s_{j} + 3n} - 15\cdot 2^{3n} +7\cdot2^{2n+1}  & \text{if   }  i=4 
    \end{cases}    
   \end{equation}

\begin{example}
We have for $ n=3,\;s_{1}=s_{2}=s_{3} =3.$\\
     \begin{equation*}
       \Gamma_{i}^{\left[3\atop{3 \atop3}\right]\times 4} =  \begin{cases}
1 & \text{if  } i = 0,        \\
 21 & \text{if   } i=1,\\
 378 & \text{if   }  i = 2,  \\
  6384 & \text{if   }  i = 3,  \\
255360 & \text{if   }  i=4 
    \end{cases}    
   \end{equation*}
See (4)
\end{example}

\begin{example}
We have for $ n=4,\;s_{1}=s_{2}=s_{3} =4.$\\
     \begin{equation*}
       \Gamma_{i}^{\left[4\atop{4 \atop{4 \atop 4}}\right]\times 4} =  \begin{cases}
1 & \text{if  } i = 0,        \\
 45 & \text{if   } i=1,\\
 1650 & \text{if   }  i = 2,  \\
  56160 & \text{if   }  i = 3,  \\
268377600 & \text{if   }  i=4 
    \end{cases}    
   \end{equation*}
   Hence  the number  $ R_{4} $ of solutions \\
  
 $(Y_1,U_{1}^{(1)},U_{2}^{(1)},U_{3}^{(1)} ,U_{4}^{(1)}, Y_2,U_{1}^{(2)},U_{2}^{(2)}, 
U_{3}^{(2)},U_{4}^{(2)},Y_3,U_{1}^{(3)},U_{2}^{(3)},U_{3}^{(3)} ,U_{4}^{(3)},Y_4,U_{1}^{(4)},U_{2}^{(4)}, U_{3}^{(4)},U_{4}^{(4)}   ) \in (\mathbb{F}_{2}[T])^{20} $ \vspace{0.5 cm}\\
 of the polynomial equations  \vspace{0.5 cm}
  \[\left\{\begin{array}{c}
 Y_{1}U_{1}^{(1)} + Y_{2}U_{1}^{(2)} + Y_{3}U_{1}^{(3)}  + Y_{4}U_{1}^{(4)} = 0  \\
    Y_{1}U_{2}^{(1)} + Y_{2}U_{2}^{(2)} +  Y_{3}U_{2}^{(3)}  + Y_{4}U_{2}^{(4)} = 0\\
      Y_{1}U_{3}^{(1)} + Y_{2}U_{3}^{(2)} +  Y_{3}U_{3}^{(3)}  + Y_{4}U_{3}^{(4)} = 0\\
Y_{1}U_{4}^{(1)} + Y_{2}U_{4}^{(2)} + Y_{3}U_{4}^{(3)}  + Y_{4}U_{4}^{(4)} = 0 
 \end{array}\right.\]
    satisfying the degree conditions \\
                   $$  degY_i \leq 3 ,
                   \quad degU_{j}^{(i)} \leq 3, \quad  for \quad 1\leq j\leq 4,  \quad 1\leq i \leq 4. $$ \\
   is equal to 

  \begin{equation*}
  2^{52}\cdot\sum_{i = 0}^{4}
  \Gamma_{i}^{\left[4\atop{4\atop {4 \atop 4}}\right]\times 4} 2^{-4i} = 2^{45}\cdot 527243
 \end{equation*}

 \end{example}

  \section{Computation of the number  $ \Gamma_{i}^{\left[s_{j}\atop{\vdots \atop s_{j}}\right]\times k} $   of   n-times persymmetric 
   $ (\sum_{j=1}^{n}s_{j}) \times k $  rank i matrices  where $s_{j} \geqslant k-1 $ for  $1\leqslant j \leqslant n $ }
  \label{sec 2} 
 
  \textbf{Conjecture :}
  Let  $s_{j} \geqslant k-1 $ for  $1\leqslant j \leqslant n .$ \\
      $ Set\quad
  (t_{1},t_{2},\ldots,t_{n} ) $ \\
 $  =  \big( \sum_{i\geq 1}\alpha _{i}^{(1)}T^{-i}, \sum_{i \geq 1}\alpha  _{i}^{(2)}T^{-i},\sum_{i \geq 1}\alpha _{i}^{(3)}T^{-i},\ldots,\sum_{i \geq 1}\alpha _{i}^{(n)}T^{-i}   \big) \in  \mathbb{P}^{n}. $\\
We state that  the number $ \Gamma_{i}^{\left[s_{j}\atop{\vdots \atop s_{j}}\right]\times k} $ of rank i n- times persymmetric matrices over the finite field  $\mathbb{F}_{2} $  of the below form
 denoted by  $D^{\left[s_{1}\atop{\vdots \atop s_{n}}\right]\times k}(t_{1},t_{2},\ldots,t_{n} ) $
   \begin{displaymath}
   \left (  \begin{array} {ccccc}
\alpha  _{1}^{(1)} & \alpha  _{2}^{(1)} & \ldots  &   \alpha_{k-1}^{(1)} &   \alpha_{k}^{(1)}  \\
\alpha  _{2}^{(1)} & \alpha  _{3}^{(1)}  & \ldots   &   \alpha_{k}^{(1)} &   \alpha_{k+1}^{(1)} \\ 
\vdots & \vdots & \vdots  & \vdots & \vdots \\
\alpha  _{s_{1}}^{(1)} & \alpha  _{s_{1}+1}^{(1)}  & \ldots  &   \alpha_{s_{1}+k-2}^{(1)} &   \alpha_{s_{1}+k-1}^{(1)}  \\ 
\hline \\
\alpha  _{1}^{(2)} & \alpha  _{2}^{(2)}  & \ldots  &   \alpha_{k-1}^{(2)} &   \alpha_{k}^{(2)}  \\
\alpha  _{2}^{(2)} & \alpha  _{3}^{(2)}  & \ldots  &   \alpha_{k}^{(2)} &   \alpha_{k+1}^{(2)} \\ 
\vdots & \vdots & \vdots & \vdots & \vdots  \\
\alpha  _{s_{2}}^{(2)} & \alpha  _{s_{2}+1}^{(2)}  & \ldots  &   \alpha_{s_{2}+k-2}^{(2)} &   \alpha_{s_{2}+k-1}^{(2)}  \\ 
\hline\\
\alpha  _{1}^{(3)} & \alpha  _{2}^{(3)}  & \ldots  &   \alpha_{k-1}^{(3)}  &   \alpha_{k}^{(3)} \\
\alpha  _{2}^{(3)} & \alpha  _{3}^{(3)}  & \ldots  &   \alpha_{k}^{(3)}&   \alpha_{k+1}^{(3)}  \\ 
\vdots & \vdots & \vdots & \vdots & \vdots \\
\alpha  _{s_{3}}^{(3)} & \alpha  _{s_{3}+1}^{(3)}  & \ldots  &   \alpha_{s_{3}+k-2}^{(3)}&   \alpha_{s_{3}+k-1}^{(3)}   \\ 
\hline \\
\vdots & \vdots & \vdots  & \vdots & \vdots \\
\hline \\
\alpha  _{1}^{(n)} & \alpha  _{2}^{(n)}  & \ldots  &   \alpha_{k-1}^{(n)} &   \alpha_{k}^{(n)}  \\
\alpha  _{2}^{(n)} & \alpha  _{3}^{(n)}  & \ldots  &   \alpha_{k}^{(n)}&   \alpha_{k+1}^{(n)}  \\ 
\vdots & \vdots & \vdots  & \vdots & \vdots \\
\alpha  _{s_{n}}^{(n)} & \alpha  _{s_{n}+1}^{(n)} & \ldots   &   \alpha_{s_{n}+k-2}^{(n)} &   \alpha_{s_{n}+k-1}^{(n)}  \\ 
\end{array} \right )  
\end{displaymath} 
  is equal to
   \begin{equation}
   \label{eq 2.1}
      \begin{cases}
1 & \text{if  } i = 0,        \\
   \Gamma_{i}^{\left[i\atop{\vdots \atop i}\right]\times (i+1)} = (2^{i+1}-1)\cdot2^{in}-3\cdot(2^{i}-1)\cdot2^{(i-1)n} + (2^{i-1}-1)\cdot2^{(i-2)n  +1}\\
   = (2^{n}-1)\cdot(2^{n+1}-1)\cdot2^{i(n+1)-2n}- (2^{n}-1)\cdot(2^{n-1}-1)\cdot2^{in-2n+1}& \text{if  } 1\leqslant i \leqslant k-1\\
2^{\sum_{j=1}^{n}s_{j}+(k-1)n} -  (2^{k}-1)\cdot2^{(k-1)n} + (2^{k-1}-1)\cdot2^{(k-2)n  +1} & \text{if  }  i=k
 \end{cases}   
   \end{equation}
where $ \Gamma_{i}^{\left[i\atop{\vdots \atop i}\right]\times (i+1)}$ denote the number of rank i n-times persymmetric matrices over $\mathbb{F}_{2}$ of the below form :  \\
  \begin{displaymath}
   \left (  \begin{array} {ccccc}
\alpha  _{1}^{(1)} & \alpha  _{2}^{(1)} & \ldots  &   \alpha_{i}^{(1)} &   \alpha_{i+1}^{(1)}  \\
\alpha  _{2}^{(1)} & \alpha  _{3}^{(1)}  & \ldots   &   \alpha_{i+1}^{(1)} &   \alpha_{i+2}^{(1)} \\ 
\vdots & \vdots & \vdots  & \vdots & \vdots \\
\alpha  _{i}^{(1)} & \alpha  _{i+1}^{(1)}  & \ldots  &   \alpha_{2i-2}^{(1)} &   \alpha_{2i-1}^{(1)}  \\ 
\hline \\
\alpha  _{1}^{(2)} & \alpha  _{2}^{(2)}  & \ldots  &   \alpha_{i}^{(2)} &   \alpha_{i+1}^{(2)}  \\
\alpha  _{2}^{(2)} & \alpha  _{3}^{(2)}  & \ldots  &   \alpha_{i+1}^{(2)} &   \alpha_{i+2}^{(2)} \\ 
\vdots & \vdots & \vdots & \vdots & \vdots  \\
\alpha  _{i}^{(2)} & \alpha  _{i+1}^{(2)}  & \ldots  &   \alpha_{2i-2}^{(2)} &   \alpha_{2i-1}^{(2)}  \\ 
\hline\\
\alpha  _{1}^{(3)} & \alpha  _{2}^{(3)}  & \ldots  &   \alpha_{i}^{(3)}  &   \alpha_{i+1}^{(3)} \\
\alpha  _{2}^{(3)} & \alpha  _{3}^{(3)}  & \ldots  &   \alpha_{i+1}^{(3)}&   \alpha_{i+2}^{(3)}  \\ 
\vdots & \vdots & \vdots & \vdots & \vdots \\
\alpha  _{i}^{(3)} & \alpha  _{i+1}^{(3)}  & \ldots  &   \alpha_{2i-2}^{(3)}&   \alpha_{2i-1}^{(3)}   \\ 
\hline \\
\vdots & \vdots & \vdots  & \vdots & \vdots \\
\hline \\
\alpha  _{1}^{(n)} & \alpha  _{2}^{(n)}  & \ldots  &   \alpha_{i}^{(n)} &   \alpha_{i+1}^{(n)}  \\
\alpha  _{2}^{(n)} & \alpha  _{3}^{(n)}  & \ldots  &   \alpha_{i+1}^{(n)}&   \alpha_{i+2}^{(n)}  \\ 
\vdots & \vdots & \vdots  & \vdots & \vdots \\
\alpha  _{i}^{(n)} & \alpha  _{i+1}^{(n)} & \ldots   &   \alpha_{2i-2}^{(n)} &   \alpha_{2i-1}^{(n)}  \\ 
\end{array} \right )  
\end{displaymath}

 \textbf{Application :}
 
   The number denoted by $ R_{q} $ of solutions \\
  
 $(Y_1,U_{1}^{(1)},U_{2}^{(1)}, \ldots,U_{n}^{(1)}, Y_2,U_{1}^{(2)},U_{2}^{(2)}, 
\ldots,U_{n}^{(2)},\ldots  Y_q,U_{1}^{(q)},U_{2}^{(q)}, \ldots,U_{n}^{(q)}   ) \in (\mathbb{F}_{2}[T])^{(n+1)q} $ \vspace{0.5 cm}\\
 of the polynomial equations  \vspace{0.5 cm}
  \[\left\{\begin{array}{c}
 Y_{1}U_{1}^{(1)} + Y_{2}U_{1}^{(2)} + \ldots  + Y_{q}U_{1}^{(q)} = 0  \\
    Y_{1}U_{2}^{(1)} + Y_{2}U_{2}^{(2)} + \ldots  + Y_{q}U_{2}^{(q)} = 0\\
    \vdots \\
   Y_{1}U_{n}^{(1)} + Y_{2}U_{n}^{(2)} + \ldots  + Y_{q}U_{n}^{(q)} = 0 
 \end{array}\right.\]
    satisfying the degree conditions \\
                   $$  degY_i \leq k-1 ,
                   \quad degU_{j}^{(i)} \leq s_{j}-1\; \text{where} \; s_{j}\geqslant k-1  \quad  for \quad 1\leq j\leq n , \quad 1\leq i \leq q $$ \\
  is equal to 
  \begin{equation}
  \label{eq 2.2}
 2^{q(\sum_{j=1}^{n}s_{j}+k) -\sum_{j=1}^{n}s_{j}-(k-1)\cdot n}\sum_{i = 0}^{k}
  \Gamma_{i}^{\left[s_{1}\atop{\vdots \atop s_{n}}\right]\times k} 2^{-iq} =  R_{q} 
 \end{equation}
\textbf{Conditional proof :}

Inspired by the results in subsection \ref{subsec  3} we proceed as follows :\\
 We assume  that $ \Gamma_{i}^{\left[s_{1}\atop{\vdots \atop s_{n}}\right]\times k} $ is equal to
   \begin{equation}
         \label{eq 2.3}
  \begin{cases}
1 & \text{if  } i = 0,        \\
 \Gamma_{1}^{\left[1\atop{\vdots \atop 1}\right]\times 2}  =   (2^{n}-1)\cdot 3 & \text{if   } i=1,\\
  \Gamma_{2}^{\left[2\atop{\vdots \atop 2}\right]\times 3} =    7\cdot 2^{2n} -9\cdot2^{n}+2 & \text{if   }  i = 2,  \\
     \Gamma_{i}^{\left[i\atop{\vdots \atop i}\right]\times (i+1)} = (2^{i+1}-1)\cdot2^{in}-3\cdot(2^{i}-1)\cdot2^{(i-1)n} +(2^{i-1}-1)\cdot2^{(i-2)n +1} \\ & \text{if  } 1\leqslant i \leqslant k-2
 \end{cases}    
   \end{equation}
 To justify our assumption \eqref{eq 2.3}  we remark that  our supposition is valid for n equal to one, two and three. \\
 \begin{enumerate}
\item \textbf{The case n =1}\\
 The number  of rank i persymmetric  matrices over $\mathbb{F}_{2}$ of the form 
    \begin{displaymath}
   \left (  \begin{array} {ccccc}
\alpha  _{1}^{(1)} & \alpha  _{2}^{(1)} & \ldots  &   \alpha_{k-1}^{(1)} &   \alpha_{k}^{(1)}  \\
\alpha  _{2}^{(1)} & \alpha  _{3}^{(1)}  & \ldots   &   \alpha_{k}^{(1)} &   \alpha_{k+1}^{(1)} \\ 
\vdots & \vdots & \vdots  & \vdots & \vdots \\
\alpha  _{s_{1}}^{(1)} & \alpha  _{s_{1}+1}^{(1)}  & \ldots  &   \alpha_{s_{1}+k-2}^{(1)} &   \alpha_{s_{1}+k-1}^{(1)}  
\end{array} \right )  
\end{displaymath} 

 is equal to \\
 
 $ (2^{i+1}-1)\cdot2^{i}-3\cdot(2^{i}-1)\cdot2^{i-1} +(2^{i-1}-1)\cdot2^{i-1} = 3\cdot2^{2i-2} $   [see (1), (2)]  \\
 
 \item \textbf{The case n =2}\\
 The number  of rank i double  persymmetric  matrices over $\mathbb{F}_{2}$ of the form 
 \begin{displaymath}
   \left (  \begin{array} {ccccc}
\alpha  _{1}^{(1)} & \alpha  _{2}^{(1)} & \ldots  &   \alpha_{k-1}^{(1)} &   \alpha_{k}^{(1)}  \\
\alpha  _{2}^{(1)} & \alpha  _{3}^{(1)}  & \ldots   &   \alpha_{k}^{(1)} &   \alpha_{k+1}^{(1)} \\ 
\vdots & \vdots & \vdots  & \vdots & \vdots \\
\alpha  _{s_{1}}^{(1)} & \alpha  _{s_{1}+1}^{(1)}  & \ldots  &   \alpha_{s_{1}+k-2}^{(1)} &   \alpha_{s_{1}+k-1}^{(1)}  \\ 
\hline \\
\alpha  _{1}^{(2)} & \alpha  _{2}^{(2)}  & \ldots  &   \alpha_{k-1}^{(2)} &   \alpha_{k}^{(2)}  \\
\alpha  _{2}^{(2)} & \alpha  _{3}^{(2)}  & \ldots  &   \alpha_{k}^{(2)} &   \alpha_{k+1}^{(2)} \\ 
\vdots & \vdots & \vdots & \vdots & \vdots  \\
\alpha  _{s_{2}}^{(2)} & \alpha  _{s_{2}+1}^{(2)}  & \ldots  &   \alpha_{s_{2}+k-2}^{(2)} &   \alpha_{s_{2}+k-1}^{(2)}  \\ 
 \end{array} \right )  
\end{displaymath} 

 is equal to \\

$ (2^{i+1}-1)\cdot2^{2i}-3\cdot(2^{i}-1)\cdot2^{2(i-1)} +(2^{i-1}-1)\cdot2^{2(i-2) +1}\\
 = 21\cdot2^{3i-4}-3\cdot2^{2i-3}$  [see(3)]
  \item \textbf{The case n =3}\\
 The number  of rank i triple persymmetric  matrices over $\mathbb{F}_{2}$ of the form 
 \begin{displaymath}
   \left (  \begin{array} {ccccc}
\alpha  _{1}^{(1)} & \alpha  _{2}^{(1)} & \ldots  &   \alpha_{k-1}^{(1)} &   \alpha_{k}^{(1)}  \\
\alpha  _{2}^{(1)} & \alpha  _{3}^{(1)}  & \ldots   &   \alpha_{k}^{(1)} &   \alpha_{k+1}^{(1)} \\ 
\vdots & \vdots & \vdots  & \vdots & \vdots \\
\alpha  _{s_{1}}^{(1)} & \alpha  _{s_{1}+1}^{(1)}  & \ldots  &   \alpha_{s_{1}+k-2}^{(1)} &   \alpha_{s_{1}+k-1}^{(1)}  \\ 
\hline \\
\alpha  _{1}^{(2)} & \alpha  _{2}^{(2)}  & \ldots  &   \alpha_{k-1}^{(2)} &   \alpha_{k}^{(2)}  \\
\alpha  _{2}^{(2)} & \alpha  _{3}^{(2)}  & \ldots  &   \alpha_{k}^{(2)} &   \alpha_{k+1}^{(2)} \\ 
\vdots & \vdots & \vdots & \vdots & \vdots  \\
\alpha  _{s_{2}}^{(2)} & \alpha  _{s_{2}+1}^{(2)}  & \ldots  &   \alpha_{s_{2}+k-2}^{(2)} &   \alpha_{s_{2}+k-1}^{(2)}  \\ 
\hline\\
\alpha  _{1}^{(3)} & \alpha  _{2}^{(3)}  & \ldots  &   \alpha_{k-1}^{(3)}  &   \alpha_{k}^{(3)} \\
\alpha  _{2}^{(3)} & \alpha  _{3}^{(3)}  & \ldots  &   \alpha_{k}^{(3)}&   \alpha_{k+1}^{(3)}  \\ 
\vdots & \vdots & \vdots & \vdots & \vdots \\
\alpha  _{s_{3}}^{(3)} & \alpha  _{s_{3}+1}^{(3)}  & \ldots  &   \alpha_{s_{3}+k-2}^{(3)}&   \alpha_{s_{3}+k-1}^{(3)}   \\ 
 \end{array} \right )  
\end{displaymath} 

 is equal to \\
$ (2^{i+1}-1)\cdot2^{3i}-3\cdot(2^{i}-1)\cdot2^{3(i-1)} +(2^{i-1}-1)\cdot2^{3(i-2) +1}\\
 = 105\cdot2^{4i-6}-21\cdot2^{3i-5} $[see (4),(5)]
\end{enumerate}

   Let $ \displaystyle  f (t_{1},t_{2},\ldots,t_{n} ) $  be the exponential sum  in $ \mathbb{P}^{n} $ defined by\\
    $(t_{1},t_{2},\ldots,t_{n} ) \displaystyle\in \mathbb{P}^{n}\longrightarrow \\
    \sum_{deg Y\leq k-1}\sum_{deg U_{1}\leq  s_{1}-1}E(t_{1} YU_{1})
  \sum_{deg U_{2} \leq s_{2}-1}E(t _{2} YU_{2}) \ldots \sum_{deg U_{n} \leq s_{n}-1} E(t _{n} YU_{n}). $\vspace{0.5 cm}\\
    Then
  $$     f (t_{1},t_{2},\ldots,t_{n} ) =
  2^{\sum_{j=1}^{n}s_{j}+k- rank\big[ D^{\left[s_{1}\atop{\vdots \atop s_{n}}\right]\times k}(t_{1},t_{2},\ldots,t_{n} )\big] } $$
    Hence  the number  $ R_{q} $    is equal to the following integral over the unit interval in $ \mathbb{K}^{n} $
    $$ \int_{\mathbb{P}^{n}} f^{q}(t_{1},t_{2},\ldots,t_{n}) dt_{1}dt _{2}\ldots dt _{n}. $$
  Observing that $ f (t_{1},t_{2},\ldots,t_{n} ) $ is constant on cosets of $ \prod_{j=1}^{n}\mathbb{P}_{s_{j}+k-1} $ in $ \mathbb{P}^{n} $\;
  the above integral is equal to 
  
  \begin{equation}
  \label{eq 2.4}
 2^{q(\sum_{j=1}^{n}s_{j}+k) -\sum_{j=1}^{n}s_{j}-(k-1)\cdot n}\sum_{i = 0}^{k}
  \Gamma_{i}^{\left[s_{1}\atop{\vdots \atop s_{n}}\right]\times k} 2^{-iq} =  R_{q} 
 \end{equation}

 From \eqref{eq 2.4} we obtain for q = 1\\
 
   \begin{equation}
  \label{eq 2.5}
 2^{k-(k-1)n}\sum_{i = 0}^{k}
 \Gamma_{i}^{\left[s_{1}\atop{\vdots \atop s_{n}}\right]\times k} 2^{-i} = 2^{\sum_{j=1}^{n}s_{j}}+2^k-1
 \end{equation}
We have obviously \\

   \begin{equation}
  \label{eq 2.6}
 \sum_{i = 0}^{k}
 \Gamma_{i}^{\left[s_{1}\atop{\vdots \atop s_{1}}\right]\times k}  = 2^{\sum_{j=1}^{n}s_{j}+(k-1)n}
 \end{equation}

From \eqref{eq 2.6} and our assumption \eqref{eq 2.3} we get :\\

\begin{eqnarray}
  \label{eqn 2.7}
 \Gamma_{k-1}^{\left[s_{1}\atop{\vdots \atop s_{1}}\right]\times k} +   \Gamma_{k}^{\left[s_{1}\atop{\vdots \atop s_{1}}\right]\times k}  =  2^{\sum_{j=1}^{n}s_{j}+(k-1)n} -\sum_{i=0}^{k-2}  \Gamma_{i}^{\left[s_{1}\atop{\vdots \atop s_{1}}\right]\times k}\\
  = 2^{\sum_{j=1}^{n}s_{j}+(k-1)n} - [2^{(k-2)n}\cdot(2^{k-1}-1) + 2^{(k-3)n}\cdot(2-2^{k-1})] \nonumber
\end{eqnarray}

From \eqref{eq 2.5} and our assumption \eqref{eq 2.3} we obtain :\\

  \begin{align}
  \label{eqn 2.8}
        2\cdot \Gamma_{k-1}^{\left[s_{1}\atop{\vdots \atop s_{1}}\right]\times k} +   \Gamma_{k}^{\left[s_{1}\atop{\vdots \atop s_{1}}\right]\times k}  =  2^{\sum_{j=1}^{n}s_{j}+(k-1)n} +2^{(k-1)n}\cdot[2^{k}-1]  -\sum_{i=0}^{k-2} 2^{k-i} \Gamma_{i}^{\left[s_{1}\atop{\vdots \atop s_{1}}\right]\times k}\\
 = 2^{\sum_{j=1}^{n}s_{j}+(k-1)n}+ 2^{(k-1)n}\cdot[2^{k}-1]  +  2^{nk-3n}\cdot[2^k-2^2]  +   2^{nk-2n}\cdot[2^2-2^{k+1}]   \nonumber
\end{align}
We deduce from  \eqref{eqn 2.7} and  \eqref{eqn 2.8} : \\
\begin{equation}
\label{eq 2.9}
 \Gamma_{k-1}^{\left[s_{1}\atop{\vdots \atop s_{1}}\right]\times k}
  =  2^{(k-1)n}\cdot[2^{k}-1] + 3\cdot 2^{nk-2n}\cdot[1-2^{k-1}]  + 2^{nk-3n+1}\cdot[2^{k-2}- 1] 
\end{equation}

\begin{equation}
\label{eq 2.10}
 \Gamma_{k}^{\left[s_{1}\atop{\vdots \atop s_{1}}\right]\times k}
  = 2^{\sum_{j=1}^{n}s_{j}+(k-1)n} -  (2^{k}-1)\cdot2^{(k-1)n} + (2^{k-1}-1)\cdot2^{(k-2)n  +1}.
\end{equation}

\maketitle 
\newpage
\tableofcontents

 \end{document}